\documentclass[12pt]{amsart}
\usepackage{graphicx}
\usepackage{amssymb}
\usepackage{epstopdf}
\usepackage{hyperref}

\def\R{\mathbf{R}}

\def\lmt{\tilde{\lambda}}

\newtheorem{innercustomgeneric}{\customgenericname}
\providecommand{\customgenericname}{}
\newcommand{\newcustomtheorem}[2]{%
  \newenvironment{#1}[1]
  {%
   \renewcommand\customgenericname{#2}%
   \renewcommand\theinnercustomgeneric{##1}%
   \innercustomgeneric
  }
  {\endinnercustomgeneric}
}

\newcustomtheorem{customthm}{Theorem}
\newcustomtheorem{customlemma}{Lemma}
\newtheorem{conj}{Conjecture}

\begin{document}
\date{}

\title[Notes]{Equilateral Chains and Cyclic Central Configurations of the Planar 5-body Problem}
\author[1]{Yiyang Deng}
\author[2]{Marshall Hampton}
\email{mhampton@d.umn.edu, Orcid ID: 0000-0002-9725-3316}

\address[1]{College of Mathematics and Statistics, Chongqing Technology and Business University, Chongqing 400067, China}
\address[2]{Department of Mathematics and Statistics, University of Minnesota Duluth, Duluth, MN 55812}

\maketitle

\begin{abstract}
Central configurations and relative equilibria are an important facet of the study of the $N$-body problem, but become very difficult to rigorously analyze for $N>3$.  In this paper we focus on a particular but interesting class of configurations of the 5-body problem: the equilateral pentagonal configurations, which have a cycle of five equal edges.  We prove a variety of results concerning central configurations with this property, including a computer-assisted proof of the finiteness of such configurations for any positive five masses with a range of rational-exponent homogeneous potentials (including the Newtonian case and the point-vortex model), some constraints on their shapes, and we determine some exact solutions for particular N-body potentials.
\end{abstract}

\section{Introduction}

In this work we consider some particular classes of relative equilibria of a planar $N$-body problem, in which $N$ point particles with non-negative masses $m_i$ interact through a central potential $U$:

$$m_i \ddot{q}_{i;j} = \frac{\partial U}{\partial q_{i;j}},  \ \ i \in \{0,\ldots N-1\}, $$

$$ U = \sum_{i<j} m_i m_j/ r_{i,j}^{A-2}$$

where $q_i \in \R^2$ is the position of particle $i$, and $r_{i,j} = |q_i - q_j|$ are the mutual distances between the particles.  The exponent $A$ is a real parameter in $(2,\infty)$.

The most interesting and important case is the Newtonian gravitational model with $A=3$, but we believe it can be useful to generalize the problem since many features of the relative equilibria do not strongly depend on the exponent $A$.  The potential can be extended to the case $A=2$ by using

$$ U = \sum_{i<k} m_i m_k \log(r_{i,k}) $$

which has been used in models of fluid vortex tubes \cite{HelmholtzV,kirchhoff_vorlesungen_1883,aref1992grobli}.


The relative equilibria (equilibria in a uniformly rotating reference frame) must satisfy the equations for a central configuration \cite{wintner41}, defined as configurations for which

\begin{equation} \label{def_cc}
\lambda (q_i - c) = \sum_{j = 1, j \neq i}^n \frac{m_j (q_i - q_j)}{r_{i,j}^A}
\end{equation}

The vector $c$ is the center of mass, 
$$c = \frac{1}{M} \sum_{i=1}^n m_i q_i,$$
with
$$M = \sum_{i=1}^{n} m_i$$
the total mass, which we will always assume to be nonzero (for the special case in which $M=0$, see \cite{celli_2005}).  The parameter $\lambda$ is real.  The masses $m_i$ are also assumed to be real, and we are primarily interested in positive masses.

In some earlier literature central configurations are also referred to as {\it permanent configurations} \cite{MacMillanBartky32, rayl39, Brumberg57}.  The study of central configurations and relative equilibria provides an avenue for progress into the N-body problem, which otherwise presents formidable difficulty.  There is a rich literature on these configurations, starting with Euler \cite{Euler_col} and Lagrange \cite{Lagrange_3} who completely characterized the relative equilibria for the Newtonian three-body problem.  The collinear three-body configurations studied by Euler were further elucidated by Moulton \cite{moulton_straight_1910}, who showed there is a unique (up to scaling) central configuration for any ordering of $N$ positive masses on a line. 

For $N>3$ it is much harder to characterize the central configurations.  One of the most basic questions is whether there are finitely many equivalence classes of central configurations for each choice of positive masses.  This question has been highlighted in the planar case by several authors \cite{chazy18, wintner41, smale98}.  It has been resolved for the Newtonian four-body problem \cite{hampton_finiteness_2005}, the four-vortex problem \cite{hampton_finiteness_2009}, and partially for the Newtonian five-body problem \cite{albouykaloshin}.  

There is much work on other interesting questions on central configurations, such as their stability.  Rather than attempt to summarize this work we recommend the excellent surveys by Moeckel \cite{LMS, moeckel_central_1990}.



\section{Equations for Central Configurations and Equilateral Chains}

Choosing two indices $i$ and $j$ we can take the inner product of (\ref{def_cc}) with $q_i - q_j$ to get

$$\lambda (q_i - c) \cdot (q_i - q_j) = \sum_{k = 1, k \neq i}^n \frac{m_k (q_i - q_k) \cdot (q_i - q_j)}{r_{ik}^A}$$
and then the left-hand side becomes
$$\lambda (q_i - c) \cdot (q_i - q_j) = \frac{\lambda}{M} \sum_{k=1}^n m_k (q_i - q_j) \cdot (q_i - q_k) = \lmt \sum_{k=1}^n m_k (r_{i,j}^2 + r_{ik}^2 - r_{j,k}^2)/2$$
in which we have introduced $\lmt = \frac{\lambda}{M}$.  The inner-products on the right-hand side can be rewritten in terms of the mutual distances as well.  After putting all the terms on one side of the equation and cancelling a factor of $1/2$, we obtain for each choice of $i \neq j$ the equations

\begin{equation} \label{elem_aceqs}
\sum_{k=1}^n m_k (r_{ik}^{-A} - \lmt) (r_{i,j}^2 + r_{ik}^2 - r_{j,k}^2) = 0
\end{equation}

If we now introduce variables $S_{i,j} = r_{ik}^{-A} - \lmt$ and  $A_{i,j,k} = r_{j,k}^2 - r_{ik}^2 - r_{i,j}^2$ we obtain the compact form 

$$f_{i,j} = \sum_{k = 1}^n m_k  S_{ik} A_{i,j,k} = 0.$$
Gareth Roberts has observed that these equations follow from the developments given in \cite{albouy_probleme_1997}; they are sometimes referred to as the `asymmetric Albouy-Chenciner equations'.  

If we combine $f_{i,j}$ and $f_{ji}$ we obtain $n(n-1)/2$ equations

$$g_{i,j,} = f_{i,j} + f_{j,i} = \sum_{k = 1}^n m_k  (S_{ik} A_{i,j,k} + S_{j,k} A_{jik}) = 0.$$
These are the equations presented as the Albouy-Chenciner equations in \cite{hampton_finiteness_2005}.  These can be interpreted kinematically as the statement that there exists a $\lambda$ for which
$$\lambda r_{i,j}^2 = - (q_i - q_j) \cdot \ddot{q}_{i,j}$$
for all $i \neq j$.

By taking the wedge product instead of the inner product, we obtain a different set of equations referred to as the Laura-Andoyer equations \cite{laura1905sulle, andoyer1906main}

\begin{equation} \label{eq:alb1}
L_{i,j} := \sum_{k \neq i,j} m_k (R_{i,k} - R_{j,k}) \Delta_{i,j,k} = 0
\end{equation}
where $\Delta_{i,j,k}$ is  the oriented area $(q_i - q_j)
\wedge (q_i - q_k)$, $q_i \in \mathbb{R}^2$ and $R_{i,j} = r_{i,j}^{-A} = (|q_i - q_j|)^{-A}$.  Sometimes the $\Delta_{i,j,k}$ will be replaced by the non-negative $D_{i,j,k} = |\Delta_{i,j,k}|$ in order to make the sign of the terms in our equations more apparent.


We define an $N$-body configuration to be an {\em equilateral chain} if at least $N-1$ consecutive distances involving all of the points are equal; we will choose the particular convention with

$$ r_{1,2} = r_{2,3} = \ldots = r_{N-2,N-1} = r_{N-1,N} $$

Similarly an {\em equilateral cyclic configuration} has a $N$ equal distances in a complete cycle, and we will choose our indexing so that 

$$ r_{1,2} = r_{2,3} = \ldots  = r_{1,N} $$

The five-body cyclic configurations generalize the rhomboidal configurations of the four-body problem which have been well studied in both the Newtonian and vortex cases \cite{long2002four,perez2007convex,hampton_vort_pairs2014,leandro2019structure,oliveira2020stability}.   The five-body configurations of a rhombus with a central mass are another interesting and well-studied extension, which contain continua of central configurations if a negative central mass is allowed \cite{roberts_continuum_1999, gidea2010symmetric, albouykaloshin, cornelio2021central}.

The Laura-Andoyer equations for the equilateral pentagon case fall into two sets of five; one of these sets contains equations involving only two of the masses:

 \begin{align}
 & m_4 \Delta_{1,3,4}(R_{1,4}-R_{1,2}) + m_5 \Delta_{1,3,5} (R_{1,2} - R_{3,5}) = 0 , \label{LA2} \\
 & m_1 \Delta_{1,2,4}(R_{1,4}-R_{1,2}) + m_5 \Delta_{2,4,5} ( R_{1,2}-R_{2,5} ) = 0 ,\nonumber \\
 & m_3 \Delta_{2,3,5}(R_{3,5}-R_{1,2}) + m_4 \Delta_{2,4,5} (R_{1,2} - R_{2,4}) = 0 , \nonumber\\
 & m_1 \Delta_{1,3,5} (R_{1,3} -R_{1,2}) +  m_2 \Delta_{2,3,5} (R_{1,2} - R_{2,5}) = 0, \nonumber\\
  & m_2  \Delta_{1,2,4} ( R_{2,4}-R_{1,2}) + m_3 \Delta_{1,3,4} (R_{1,2}- R_{1,3} ) = 0,  \nonumber  
 \end{align}
 
 while the remaining equations involve three masses:
 \begin{align*}
 & m_3 \Delta_{1,2,3} (R_{1,3} -R_{1,2} )  + m_4  \Delta_{1,2,4} (R_{1,4}- R_{2,4})+ m_5  \Delta_{1,2,5} (R_{1,2} - R_{2,5} ) = 0, \\
 & m_2 \Delta_{1,2,5} (R_{2,5} -R_{1,2} )  + m_3  \Delta_{1,3,5} (R_{3,5}- R_{1,3})+ m_4  \Delta_{1,4,5} (R_{1,2} - R_{1,4} ) = 0, \\
 & m_1 \Delta_{1,2,3} (R_{1,2} - R_{1,3} )  + m_4  \Delta_{2,3,4} (R_{2,4}- R_{1,2})+ m_5  \Delta_{2,3,5} (R_{2,5} - R_{3,5} ) = 0, \\
 & m_1 \Delta_{1,3,4} (R_{1,3} -R_{1,4} )  + m_2  \Delta_{2,3,4} (R_{1,2}- R_{2,4})+ m_5  \Delta_{3,4,5} (R_{3,5} - R_{1,2} ) = 0, \\
 & m_1 \Delta_{1,4,5} (R_{1,4} -R_{1,2} )  + m_2  \Delta_{2,4,5} (R_{2,4}- R_{2,5})+ m_3  \Delta_{3,4,5} (R_{1,2} - R_{3,5} ) = 0, 
 \end{align*}
 
 If we normalize the configurations by choosing $q_1 = (-1/2, 0)$ and $q_2 = (1/2,0)$ (so $r_{1,2} = 1$), then the $\Delta_{i,j,k}$ are
 
 $$ \Delta_{1,2,3} = y_3, \ \ \ \ \Delta_{1,2,4} = y_4, \ \ \ \ \Delta_{1,2,5} = y_5$$
 
 $$ \Delta_{1,3,4} = -x_{4} y_{3} + x_{3} y_{4} + \frac{1}{2} (y_{4} - y_{3} ) $$
 
 $$ \Delta_{1,3,5} =  -x_{5} y_{3} + x_{3} y_{5} + \frac{1}{2} (y_{5} - y_{3} )$$
 
 $$ \Delta_{1,4,5} = -x_{5} y_{4} + x_{4} y_{5} + \frac{1}{2} (y_{5} - y_{4} )$$
 
 $$ \Delta_{2,3,4} = -x_{4} y_{3} + x_{3} y_{4} + \frac{1}{2} (y_{3} - y_{4} )$$
 
 $$ \Delta_{2,3,5} = -x_{5} y_{3} + x_{3} y_{5} + \frac{1}{2} (y_{3} - y_{5} )$$
 
 $$ \Delta_{2,4,5} = -x_{5} y_{4} + x_{4} y_{5} + \frac{1}{2} (y_{4} - y_{5} )$$
 
 $$ \Delta_{3,4,5} = -x_{4} y_{3} + x_{5} y_{3} + x_{3} y_{4} - x_{5} y_{4} - x_{3} y_{5} + x_{4} y_{5} $$

\section{Finiteness results}


To study the finiteness of the planar equilateral configurations we use the asymmetric Albouy-Chenciner equations $f_{i,j}=0$ and the Cayley-Menger determinants for all four-point subconfigurations

\begin{equation} \label{CM}
C_{i,j,k,l} = det \left ( \begin{array}{ccccc} 0 & 1 & 1 & 1 & 1  \\ 
1 & 0 &  r_{i,j}^2 & r_{i,k}^2 & r_{i,l}^2  \\
1 & r_{i,j}^2 & 0 & r_{j,k}^2 & r_{j,l}^2  \\
1 & r_{i,k}^2 &  r_{j,k}^2 & 0 & r_{k,l}^2 \\
1 & r_{i,l}^2 &  r_{j,l}^2 & r_{k,l}^2 & 0 \\
 \end{array} \right ) = 0
\end{equation}

(with distinct $i$, $j$, $k$, and $l$).

Following a very similar procedure to that described in \cite{hampton_finiteness_2011}, with the assistance of the software Sage \cite{sage_2020}, Singular \cite{DGPS21}, and Gfan \cite{gfan} we find relatively easily that there are finitely many equilateral central configurations for $A=2$ and $A=3$, for any positive masses.  The tropical prevariety of the system is much simpler in the vortex case $A=2$, having only 22 generating rays compared to 37 in the Newtonian case.   In fact the Newtonian case generalizes for potential exponents greater than 3 as well.

Because Gfan cannot currently compute a prevariety for variable exponents, we used a slight variation of the Albouy-Chenciner equations to compute and to analyze the initial form systems for rational exponents $A>2$.  Let $Q_{i,j} = r_{i,j}^{-A+2}$, so that

$$ f_{i,j} = \sum_{k = 1}^n m_k  S_{ik} A_{i,j,k} = \sum_{k = 1}^n m_k  (Q_{i,j}r_{i,j}^{-2} - \lmt) A_{i,j,k} $$

We used Gfan to compute the prevariety of these equations with the normalization $\lmt = 1$, and cleared denominators to obtain a polynomial system in $r_{i,j}$ and $Q_{i,j}$, combined with the four-point Cayley-Menger determinants. 

For rational $A > 2$, the rays of the tropical prevariety fall into 9 equivalence classes under the action of the cyclic group $C_5$ on the indices of the $r_{i,j}$.  Since we have chosen the equilateral configurations to have $r_{i,j} = r_{i+1, j+1}$, this action fixes our first coordinate $r_{i,j}$ and cyclically permutes the remaining five distances.  In the table below, the coordinate exponents are in the order $(r_{1,2}, r_{1,3}, r_{1,4}, r_{2,4}, r_{2,5}, r_{3,5})$.  Six of these rays are independent of $A$.

\begin{center}
\begin{tabular}{|c|c|c|} \hline
Label & Ray Representative & Multiplicity  \\ \hline
$h_1$ &  $(-1,-1,-1,-1,-1,-1)$ & 1 \\ \hline
$h_2$ & $(1,1,1,1,1,1)$ & 1 \\ \hline
$h_3$ & $(0,1,1,1,1,1)$ & 1 \\ \hline
$h_4$ & $(0,0,1,1,1,1)$ & 5 \\ \hline
$h_5$ & $(1,0,1,1,1,1)$ & 5 \\ \hline
$h_6$ & $(1, 0, 1, 0, 1, 1)$ & 5 \\ \hline
$h_7$ & $(A-2, -2, A-2, -2, A-2, A-2)$ & 5 \\ \hline
$h_8$ & $(A-2, -2, A-2, 0, A-2, A-2)$ & 10 \\ \hline
$h_9$ & $(A-2, -2, A-2, A-2, A-2, A-2)$ & 5 \\ \hline
\end{tabular}
\end{center}


Because of the balance condition for tropical varieties \cite{Maclagan}, we can restrict our analysis to cones that intersect the half-space containing exponent vectors with a non-negative sum.  This excludes the first ray in our list.  Again after reducing by the $C_5$ symmetry we have cones

\begin{center}
\begin{tabular}{|c|c|} \hline
Label & Representative Cone Rays   \\ \hline
$C_1$  &  \{(0, 0, 1, 1, 1, 1)\}  \\ \hline
$C_2$  &  \{(0, 1, 1, 1, 1, 1)\}  \\ \hline
$C_3$ &  \{(1, 0, 1, 0, 1, 1)\}  \\ \hline
$C_4$  &  \{(1, 0, 1, 1, 1, 1)\}  \\ \hline
$C_5$  &  \{(1, 1, 1, 1, 1, 1)\}  \\ \hline
$C_6$  &  \{(A-2, -2, A-2, -2, A-2, A-2)\}  \\ \hline
$C_7$  &  \{(A-2, -2, A-2, 0, A-2, A-2)\}  \\ \hline
$C_8$  &  \{(A-2, -2, A-2, A-2, 0, A-2)\}  \\ \hline
$C_9$  &  \{(A-2, -2, A-2, A-2, A-2, A-2)\}  \\ \hline
$C_{10}$  &  \{(0, 0, 1, 1, 1, 1), (0, 1, 1, 0, 1, 1)\}  \\ \hline
$C_{11}$  &  \{(0, 0, 1, 1, 1, 1), (0, 1, 1, 1, 0, 1)\}  \\ \hline
$C_{12}$  &  \{(0, 0, 1, 1, 1, 1), (0, 1, 1, 1, 1, 1)\}  \\ \hline
$C_{13}$  &  \{(0, 0, 1, 1, 1, 1), (1, 1, 1, 1, 1, 1)\}  \\ \hline
$C_{14}$  &  \{(0, 1, 1, 1, 1, 1), (1, 1, 1, 1, 1, 1)\}  \\ \hline
$C_{15}$  &  \{(1, 0, 1, 1, 1, 1), (1, 1, 1, 1, 1, 1)\}  \\ \hline
$C_{16}$  &  \{(1, 0, 1, 0, 1, 1), (A-2, -2, A-2, 0, A-2, A-2)\}  \\ \hline
$C_{17}$  &  \{(1, 0, 1, 0, 1, 1), (A-2, 0, A-2, -2, A-2, A-2)\}  \\ \hline
$C_{18}$  &  \{(1, 0, 1, 1, 1, 1), (A-2, -2, A-2, A-2, A-2, A-2)\}  \\ \hline
$C_{19}$  &  \{(A-2, -2, A-2, -2, A-2, A-2), (A-2, -2, A-2, 0, A-2, A-2)\}  \\ \hline
$C_{20}$  &  \{(A-2, -2, A-2, -2, A-2, A-2), (A-2, 0, A-2, -2, A-2, A-2)\}  \\ \hline
$C_{21}$  &  \{(0, 0, 1, 1, 1, 1), (0, 1, 1, 1, 1, 1), (1, 1, 1, 1, 1, 1)\}  \\ \hline
$C_{22}$  &  \{(1, 0, 1, 0, 1, 1), (A-2, -2, A-2, -2, A-2, A-2),  \\ 
      & (A-2, -2, A-2, 0, A-2, A-2), (A-2, 0, A-2, -2, A-2, A-2)\} \\ \hline
\end{tabular}
\end{center}

 Remarkably, for each cone in the tropical prevariety the initial form polynomials factor enough that we can compute the elimination ideal of the initial form system in the ring $\mathbb{Q} [m_1, m_2, m_3, m_4, m_5]$, using Singular \cite{DGPS420} within Sage \cite{sagemath91}, without the need to specialize to a particular value of $A$ (apart from the condition that $A \neq 2$. To rule out a nontrivial (i.e. nonmonomial) initial form ideal we only need to assume that no subset of the masses has a vanishing sum, including the total mass.  

It seems possible that this formulation of the equations with the $Q_{i,j}$ may be useful in studying central configurations in other contexts.

We can include the case of rational $A$ in this result, since if $A = p/q$ then we can use the polynomial conditions $Q_{i,j}^q r_{i,j}^p - r_{i,j}^{2 q}=0$ to define the $Q_{i,j}$.

Altogether this gives a computer-assisted proof of

\begin{customthm}1 For nonzero masses with nonzero subset sums (i.e. $m_i + m_j \neq 0$, $m_i + m_j + m_k \neq 0$, $m_i + m_j + m_k + m_l \neq 0$, $m_1 + m_2 + m_3 + m_4 + m_5  \neq 0$), there are finitely many planar equilateral 5-body central configurations for any rational potential exponent $A \ge 2$.
\end{customthm}

This result strongly suggests the conjecture that there are finitely many planar equilateral 5-body central configurations for any real $A \ge 2$, but the proof of that would require different methods.

\section{The symmetric case}

In this section we impose the further restrictions of an axis of symmetry $r_{1,3} = r_{2,5}$, $r_{1,5} = r_{2,3}$, and $r_{1,4} = r_{2,4}$. 

In addition to normalizing the size of the configuration with $r_{1,2} = 1$, we can choose cartesian coordinates $q_1 = (-1/2, 0)$, $q_2 = (1/2,0)$, $q_3 = (x_3, y_3)$, $q_4 = (0, y_4)$, and $q_5 = (-x_3, y_3)$.  The equilateral constraints in these coordinates become:

$$ 4 x_3^2 - 4 x_3 + 4 y_3^2 - 3 = 0 $$

$$ x_3^2 + y_3^2 - 2 y_3 y_4 + y_4^2 - 1  = 0 $$

We can parameterize the configurations by $y_4$, in terms of which

$$ y_3 = \frac{8 \, y_{4}^{3} + 2 \, y_{4} \pm \sqrt{-16 \, y_{4}^{4} + 56 \, y_{4}^{2} + 15}}{4 \, {\left(4 \, y_{4}^{2} + 1\right)}} $$

$$ x_3 = \frac{4 \, y_{4}^{2} + 1 \pm 2 \, y_4 \sqrt{-16 \, y_{4}^{4} + 56 \, y_{4}^{2} + 15} }{4 \, {\left(4 \, y_{4}^{2} + 1\right)}} $$

Note that the choices of sign must be the same, giving us two curves of configurations.  We will refer to the positive choice of sign as branch A, and the other as branch B.

The wedge products $\Delta_{i,j,k}$ in these coordinates are

\begin{align*}
 \Delta_{1,2,3} & = y_3  \\ 
 \Delta_{1,2,4} & = y_4   \\ 
 \Delta_{1,3,4} & = x_3 y_4 + (y_4 - y_3)/2   \\ 
 \Delta_{1,3,5} & = 2 y_3  x_3   \\ 
 \Delta_{1,4,5} & = x_3 y_4 - (y_4 - y_3)/2   \\ 
 \Delta_{3,4,5} & = 2 x_3 (y_4 - y_3)     \\ 
\end{align*}

\begin{figure}
\includegraphics[width=3.5in]{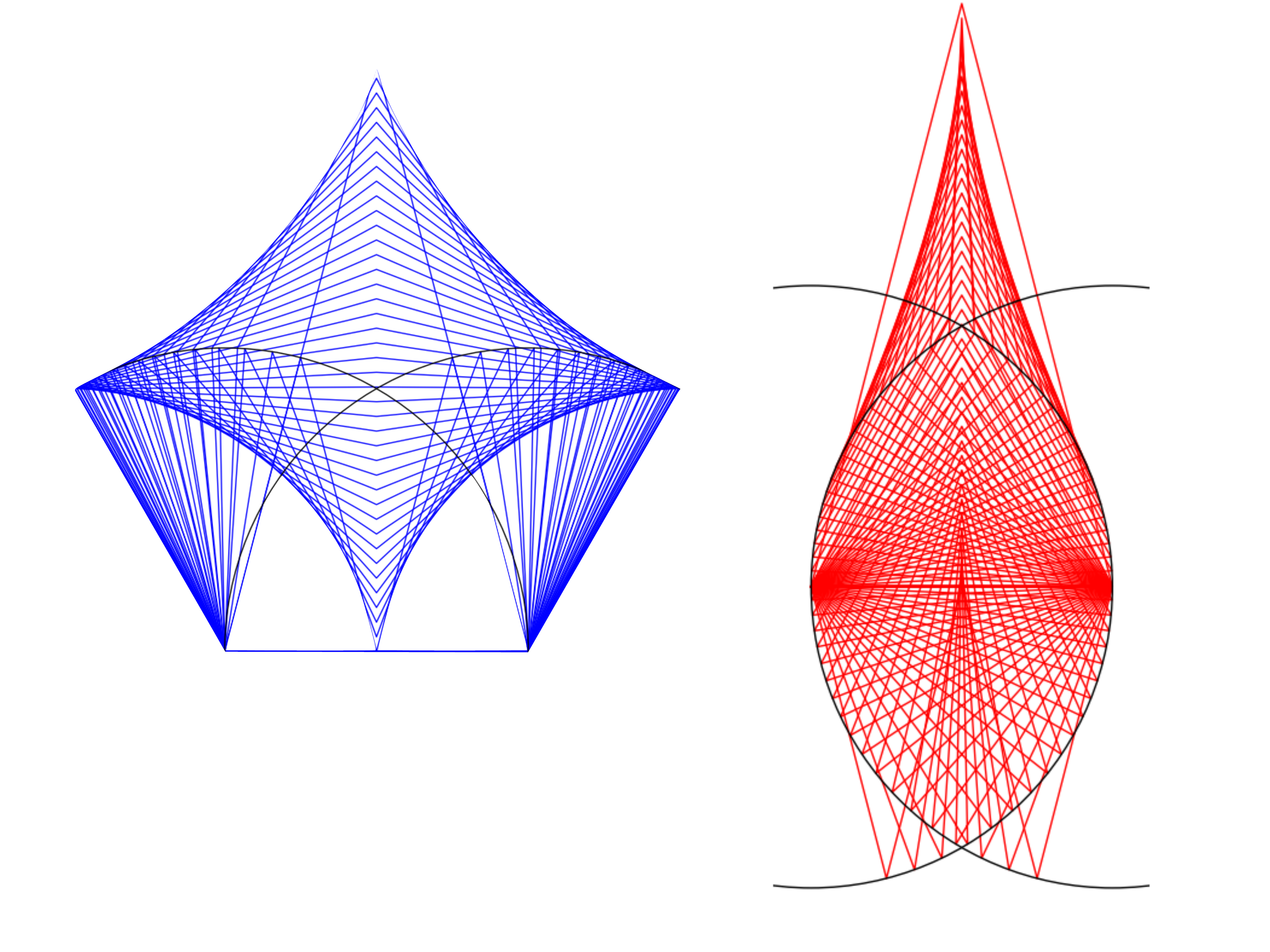}
\caption{Normalized symmetric equilateral pentagons.  Branch A is in blue, branch B in red.} \label{equifig}
\end{figure}

These configurations are shown in Figure \ref{equifig}.

For the vortex case $A=2$, it is possible to compute a Groebner basis of the system used in the finiteness proof.  This basis shows there are only two possible symmetric equilateral configurations: the regular pentagon with equal masses, and a configuration (using normalized masses $m_1 = m_2 = 1$) with $m_4$ satisfying 

$$ 64  m_4^9 - 752  m_4^8 + 2316  m_4^7 - 109  m_4^6 - 2830  m_4^5 + 45  m_4^4 + 1362  m_4^3 + 215  m_4^2 - 149  m_4 - 17 = 0 $$

which has a single positive root at $m_4 \approx 0.34199$ and then $m_3=m_5 \approx 2.32$.

We were also able to compute a Groebner basis for $A=4$, with the mass polynomial

\begin{align*} & 12288  m_4^{16} - 232064  m_4^{15} + 636883  m_4^{14} + 5616221  m_4^{13} + 2342977  m_4^{12} - 15626678  m_4^{11} \\
& - 6546497  m_4^{10}  + 17143788  m_4^9 - 1407668  m_4^8 - 5326884  m_4^7 + 456601  m_4^6 + 2374416  m_4^5 \\
& - 239673  m_4^4 - 387130  m_4^3 - 33431  m_4^2 + 25519  m_4 + 957  = 0
\end{align*}

\begin{figure}[h!b]
\includegraphics[width=3in]{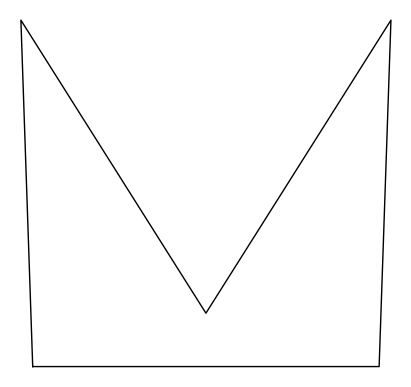}
\caption{Symmetric vortex ($A=2$) central configuration}
\end{figure}

The Laura-Andoyer equations for the symmetric case are

$$ L_{1,3} = m_3 (1 - R_{3,5}) \Delta_{1,3,5} + m_4 (R_{1,4} - 1) \Delta_{1,3,4} = 0 $$

$$ L_{1,4} = m_1 (1 - R_{1,4}) \Delta_{1,2,4} + m_3 (R_{1,3} - 1) \Delta_{1,3,4} 
= 0 $$

$$ L_{1,5} = m_1 (1 - R_{1,3}) \Delta_{1,2,3} + m_3 (R_{1,3} - R_{3,5}) \Delta_{1,3,5} + m_4 (R_{1,4} - 1)\Delta_{1,4,5}  = 0 $$

$$ L_{3,4} = m_1 \left [ (R_{1,3} - R_{1,4}) \Delta_{1,3,4} + (1 - R_{1,4}) \Delta_{2,3,4} \right ] + m_3 (R_{3,5} - 1) \Delta_{3,4,5} = 0 $$

We can highlight the linearity of these equations in the masses by forming the mass-coefficient matrix:

$$ \left(\begin{array}{ccc}
0 & ( 1- R_{35}) \Delta_{1,3,5} & (R_{1,4} - 1) \Delta_{1,3,4} \\
(1 - R_{1,4}) \Delta_{1,2,4} & (R_{1,3} - 1) \Delta_{1,3,4} & 0 \\
(1 - R_{1,3}) \Delta_{1,2,3} & (R_{1,3} - R_{3,5}) \Delta_{1,3,5} & (R_{1,4} - 1)\Delta_{1,4,5}  \\
(R_{1,3} - R_{1,4}) \Delta_{1,3,4} + (1 - R_{1,4}) \Delta_{1,4,5}  & (R_{3,5} - 1) \Delta_{3,4,5} & 0
\end{array}\right) $$

We can row-reduce this a little to get

\begin{equation} \label{MassCoef}
 \left(\begin{array}{ccc}
0 & ( 1- R_{35} )\Delta_{1,3,5} & (R_{1,4} - 1) \Delta_{1,3,4} \\
(1 - R_{1,4}) \Delta_{1,2,4} & (R_{1,3} - 1) \Delta_{1,3,4} & 0 \\
(1 - R_{1,3}) \Delta_{1,2,3} & \left [ (R_{1,3} - R_{3,5})  -  (1 - R_{3,5}) \frac{\Delta_{1,4,5}}{\Delta_{1,3,4}} \right ] \Delta_{1,3,5}& 0  \\
(R_{1,3} - R_{1,4}) \Delta_{1,3,4} + (1 - R_{1,4}) \Delta_{1,4,5}  & (R_{3,5} - 1) \Delta_{3,4,5} & 0
\end{array}\right) 
\end{equation}

This matrix must have a kernel vector of masses in the positive orthant.  This imposes many constraints.  To simplify the analysis of these constraints, we assume (without loss of generality) that $y_4 \ge 0$.  This convention means that $\Delta_{1,2,4} \ge 0$.  It is immediate from equation $L_{1,4}$ that $y_4$ cannot be zero, so we can assume that $\Delta_{1,2,4}$ is strictly positive.

\begin{figure}[h!b]
\includegraphics[width=5in]{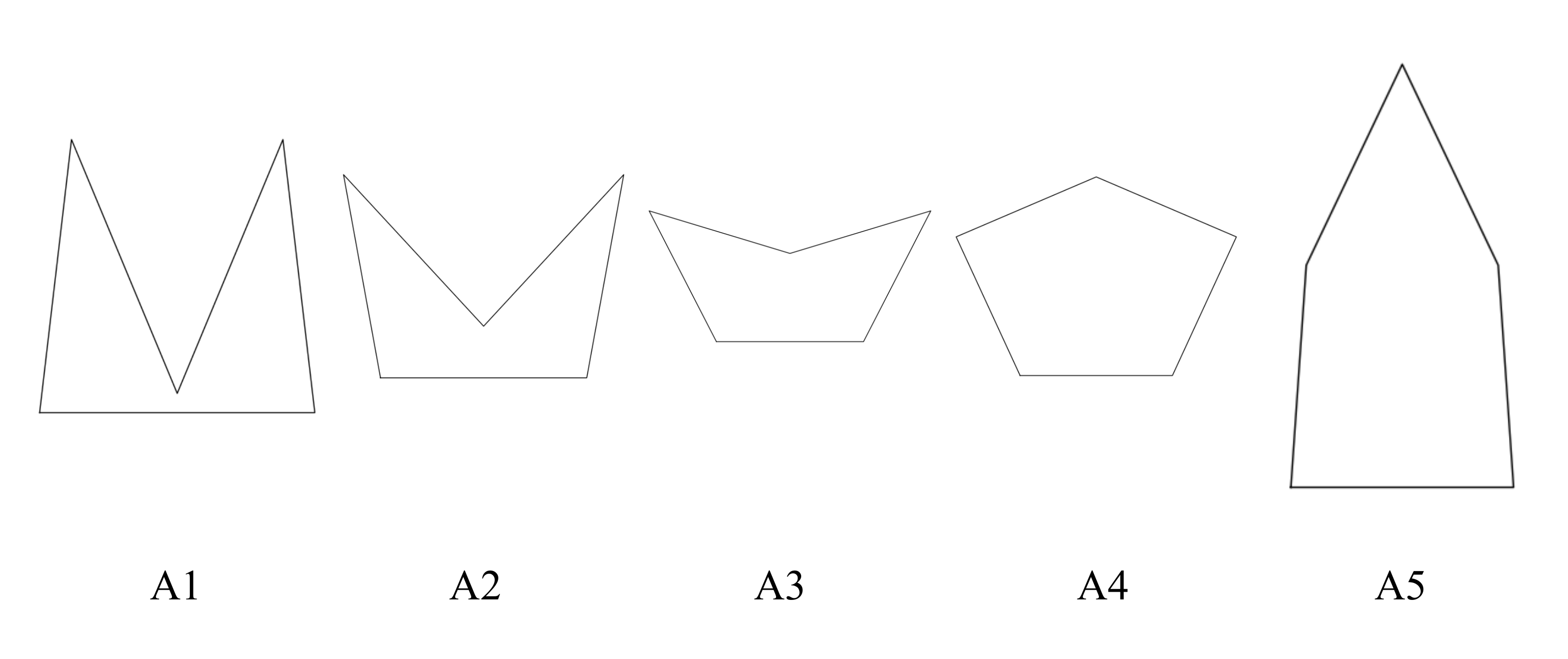}
\caption{Sign type representatives of A branch configurations}
\label{branch_As}
\end{figure}

In terms of the signs of the $\Delta_{i,j,k}$ and the magnitude  of the mutual distances relative to $r_{1,2}$, there are five cases for the branch A configurations.  Representatives of these are shown in Figure \ref{branch_As}.

All of the branch A configurations have $\Delta_{1,2,3} > 0$, $\Delta_{1,2,4} > 0 $, $\Delta_{1,3,5} > 0$, $\Delta_{1,4,5} > 0$, and $r_{1,3} > \frac{\sqrt{6}}{2} > 1$.  The distinguishing geometric properties of the branch A configurations are:

\begin{enumerate}

\item[A1)] Concave, $\Delta_{1,3,4} < 0$, $\Delta_{3,4,5} < 0$, $r_{1,4} < 1$, $r_{3,5} < 1$.

\item[A2)] Concave, $\Delta_{1,3,4} < 0$, $\Delta_{3,4,5} < 0$, $r_{1,4} < 1$, $r_{3,5} > 1$.

\item[A3)] Concave, $\Delta_{1,3,4} > 0$, $\Delta_{3,4,5} < 0$,  $r_{1,4} < 1$, $r_{3,5} > 1$.

\item[A4)] Convex, $\Delta_{1,3,4} > 0$, $\Delta_{3,4,5} > 0$,  $r_{1,4} > 1$, $r_{3,5} > 1$.

\item[A5)] Convex, $\Delta_{1,3,4} > 0$, $\Delta_{3,4,5} > 0$, $r_{1,4} > 1$, $r_{3,5} > 1$.
\end{enumerate}

\begin{customthm}{2}
The only possible branch A central configurations are of type A2 or A4.  There is a unique type A2 central configuration for all $A \ge 2$. 
\end{customthm}

\begin{proof}

For type A1 configurations, we can immediately verify that $L_{1,3} < 0$, so no such central configurations are possible.  This is also true for the borderline A1/A2 configuration with $r_{3,5} = 1$.  

The type A2 configurations contain an interesting central configuration.  

For the mass-coefficient matrix \ref{MassCoef} to have nonzero mass solutions, the following minor must vanish:
\begin{align*} 
  F(y_4) & = \det   \left(\begin{array}{cc}
(1 - R_{1,4}) \Delta_{1,2,4} & (R_{1,3} - 1) \Delta_{1,3,4} \\
(R_{1,3} - R_{1,4}) \Delta_{1,3,4} + (1 - R_{1,4}) \Delta_{1,4,5}  & (R_{3,5} - 1) \Delta_{3,4,5} 
\end{array}\right)  \\
  & = (1 - R_{1,4})(R_{3,5} - 1) \Delta_{1,2,4} \Delta_{3,4,5} \\
  & \ \ + (1 - R_{1,3}) \Delta_{1,3,4} ((R_{1,3} - R_{1,4}) \Delta_{1,3,4} + (1 - R_{1,4}) \Delta_{1,4,5}) = 0
\end{align*}

We can prove the existence of a type A2 central configuration for any $A \ge 2$ by examining the sign of $F$ at the regional endpoints.  At the lower endpoint,  where $y_4 = \frac{2 - \sqrt{3}}{2})$, the points $1$,$2$,$3$, and $5$ form a square, with $r_{3,5} = 1$, $r_{1,4} = \sqrt{ 2 - \sqrt{3}}$, $\Delta_{1,3,4} = \frac{1 - \sqrt{3}}{2} < 0$, and $\Delta_{1,4,5} = \frac{1}{2}$.    Since $R_{3,5} = 1$, the first portion of $F$ is zero, which makes it elementary to check the sign of the remaining terms and see that $F(\frac{2 - \sqrt{3}}{2}) > 0$.   

At the other endpoint, $y_4 = \frac{\sqrt{5 - 2 \sqrt{5}}}{2}$, and the points $1$,$4$, and $3$ are collinear (as are points $2$, $4$, and $5$) so $\Delta_{1,3,4} = 0$.  Only the first portion of $F$ is nonzero for this configuration, and it is straightforward to see that $F(\frac{\sqrt{5 - 2 \sqrt{5}}}{2}) < 0$.

We can prove the uniqueness of the type A2 central configuration for each $A \ge 2$ using interval arithmetic, simply evaluating $F$ and $\frac{d F}{d y_4}$ for intervals of $y_4$ and $A$.  The lack of any common zeros shows that no bifurcations occur, and it suffices to check the uniqueness for $A=2$, for which we have a Groebner basis.

For type A3 configurations, and the boundary case with $D_{1,3,4} = 0$, it is immediate that $L_{1,3} > 0$, so no such central configurations are possible.

Type A4 configurations include the regular pentagon, which is a central configuration for equal masses for all potential exponents $A$.  Remarkably, there is a bifurcation at $A_c \approx 3.12036856$.  It appears that for $A \in [2, A_c]$ the regular pentagon is the only central configuration of type A4, and for $A>A_c$ there are 3 type A4 central configurations (including the regular pentagon).  For large $A$, these two new central configurations converge to the A3/A4 boundary case for which $D_{3,4,5} = 0$, and the house-like configuration with $r_{3,5} = 1$ and $y_4 = 1 + \frac{\sqrt{3}}{2}$.   With interval arithmetic we can prove the uniqueness of the regular pentagon for $A \in [2,3]$, but we do not have an exact value for the bifurcation value $A_c$.  

For type A5 configurations once again $L_{1,3} < 0$, and no such central configurations are possible.

\end{proof}

For the branch B configurations there are three subcases of convex configurations and two concave, along with some exceptional cases on the borders between them.  For all branch B configurations $\Delta_{1,4,5} < 0$ and $r_{3,5} < 1$.  The subcases are:

\begin{enumerate}

\item[B1)] Convex, $\Delta_{1,2,3} < 0$, $\Delta_{1,3,4} > 0$, $\Delta_{1,3,5} < 0$,  $\Delta_{3,4,5} > 0$, $r_{1,3} > 1$, and $r_{1,4} < 1$. 

\item[B2)] Convex, $\Delta_{1,2,3} < 0$, $\Delta_{1,3,4} > 0$, $\Delta_{1,3,5} > 0$, $\Delta_{3,4,5} < 0$, $r_{1,3} < 1$, and $r_{1,4} < 1$.

\item[B3)] Convex, $\Delta_{1,2,3} > 0$, $\Delta_{1,3,4} < 0$, $\Delta_{1,3,5} < 0$, $\Delta_{3,4,5} < 0$, $r_{1,3} < 1$, and $r_{1,4} > 1$.

\item[B4)] Concave, $\Delta_{1,2,3} > 0$, $\Delta_{1,3,4} \geq 0$, $\Delta_{1,3,5} < 0$, $\Delta_{3,4,5} < 0$, $r_{1,3} < 1$, and $r_{1,4} > 1$.

\item[B5)] Concave, $\Delta_{1,2,3} > 0$, $\Delta_{1,3,4} > 0$, $\Delta_{1,3,5} > 0$, $\Delta_{3,4,5} > 0$, $r_{1,3} > 1$, and $r_{1,4} > 1$.
\end{enumerate}

\begin{figure}[h!b]
\includegraphics[width=5in]{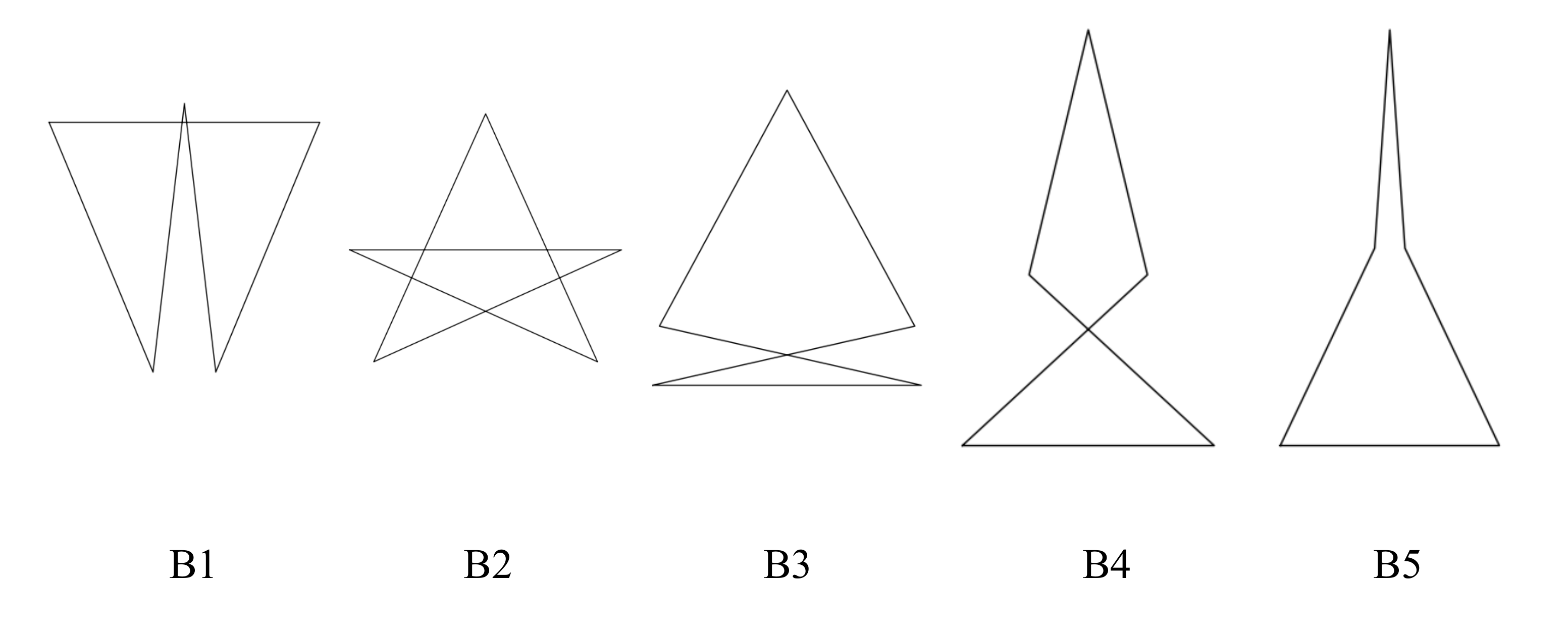}
\caption{Sign type representatives of B branch configurations}
\end{figure}

\begin{customthm}{3}
The only possible branch B central configuration is the regular pentagon for all $A \ge 2$ (type B2).
\end{customthm}

\begin{proof}

The case B1 has no solutions with positive masses, as both terms in equation $L_{1,3}$ are positive.

The case B2 contains the regular pentagon. Using interval arithmetic with the function $F(y_4)$ and its derivative we can verify that there are no other type B2 central configurations for all $A \ge 2$.

The case B3 has no solutions with positive masses, as both terms in equation $L_{1,3}$ are positive.

The case B4 has no solutions with positive masses, as both terms in equation $L_{1,4}$ are positive.

The case B5 has no solutions with positive masses, as both terms in equation $L_{1,3}$ are negative.

\end{proof}

\section{General planar equilateral pentagons}

We know that two diagonals are always greater than the four exterior edges for planar convex four-body central configurations. For strictly convex planar five-body problem, Chen and Hsiao \cite{chen2018strictly} showed that at least two interior edges are less than the exterior edges if the five bodies form a central configuration. They also showed numerical examples of strictly convex central configurations with five bodies that have either one or two interior edges less than the exterior edges. However for convex planar equilateral 5-body central configurations we have the following result:

\begin{customlemma}1  For planar convex equilateral 5-body central configurations all interior edges are greater than the exterior edges.
\end{customlemma}

\begin{proof} We consider the Laura-Andoyer equations involving only two of the masses (equations (\ref{LA2})). For the convex case we know that  $ \Delta_{1,2,3}, \Delta_{1,2,4}, \Delta_{1,2,5},  \Delta_{1,3,4}, \Delta_{1,3,5}, \Delta_{1,4,5}, \Delta_{2,3,4}, \Delta_{2,3,5}$ and  $\Delta_{3,4,5}$ are all positive.  

There is at least one interior edge greater than the exterior edges for any planar equilateral 5-body configuration. Without loss of generality, let $r_{1,4} > r_{1,2}$, and then $R_{1,4} < R_{1,2}$ and $R_{1,4} - R_{1,2} < 0$ . From the first equation of (\ref{LA2}) we must have $R_{1,2} - R_{3,5} >0$. So $R_{1,2} > R_{3,5}$ and $r_{3,5} > r_{1,2}$. Similarly, from the third equation, we have $r_{2,4} > r_{1,2}$; from the second equation above, we get $r_{1,3} > r_{1,2}$; from the fifth equation, we obtain $r_{2,5} > r_{1,2}$.
Thus, we have
$$r_{1,3}, r_{1,4}, r_{2,4}, r_{2,5}, r_{3,5} > r_{1,2}.$$
\end{proof}

The five simpler Laura-Andoyer equations (\ref{LA2}) can be used to further restrict the possible configurations of planar equilateral 5-body central configurations. 

The allowed regions fall into three classes:

\begin{enumerate}
    \item[Region I:] defined by $r_{1,2} < r_{i,j}$ and containing the regular pentagon (in 12345 order)
    
    \item[Region II:] defined by $r_{1,2} > r_{i,j}$ and containing the regular pentagon (in 13524 order)
    
    \item[Region III:] five disjoint regions, which are equivalent under permutations.  These are concave configurations.  The region with point 5 in the interior is defined by  $\theta_{1,2,3} + \theta_{2,3,4} \le 3 \pi$, $\theta_{1,2,3} \le 5 \pi/3$, $\theta_{2,3,4} \le 5 \pi/3$, $\Delta_{1,3,5} \ge 0$, and $\Delta_{2,4,5} \ge 0$.
    
\end{enumerate}

\begin{center}
\begin{figure} 
\includegraphics[width=4in]{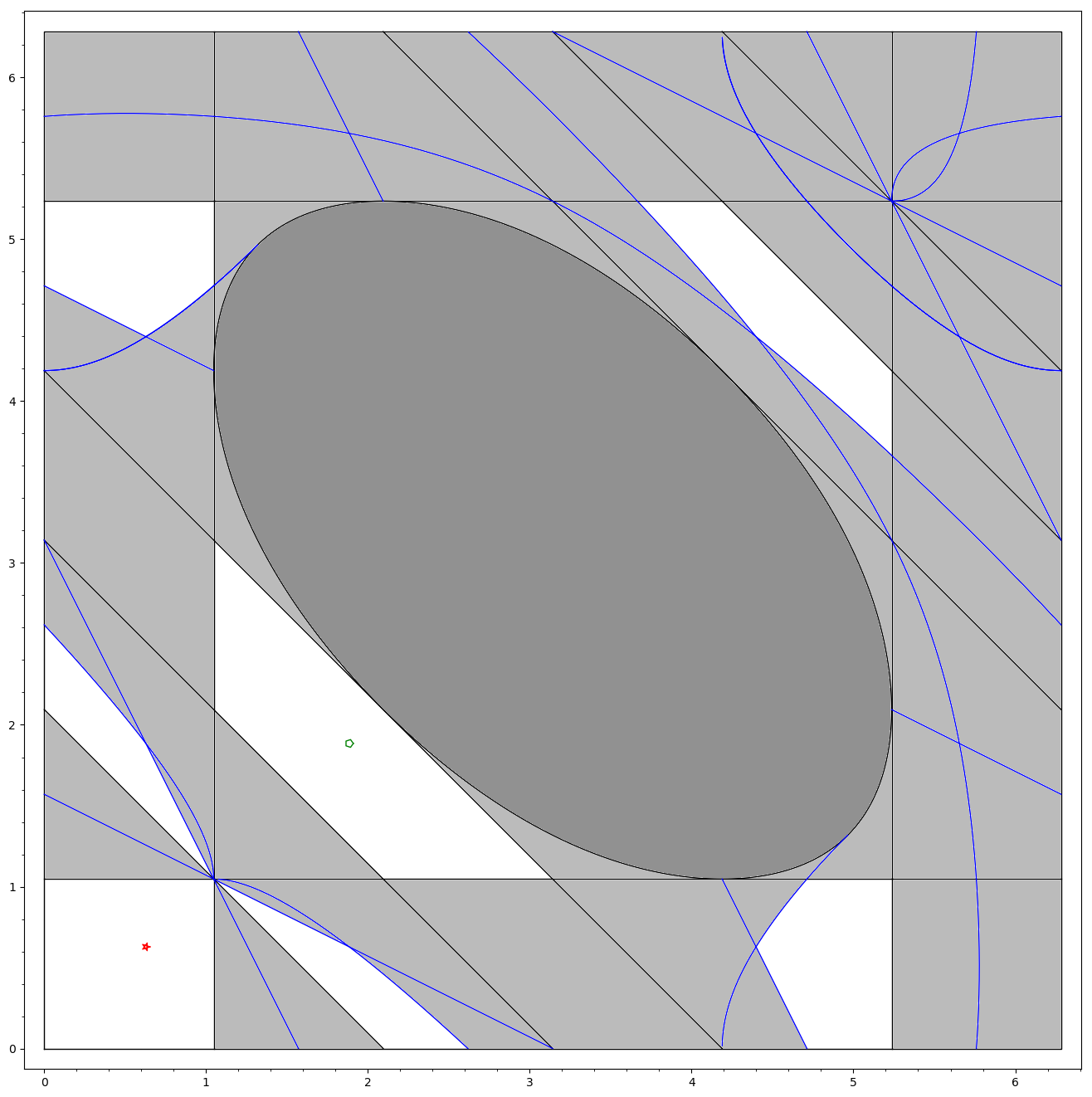}
\caption{Regions of allowable angles $\theta_{1,2}$ and $\theta_{2,3}$ for positive mass equilateral central configurations (white).  The more heavily shaded central oval consists of angles that are not geometrically realizable. The regular pentagon configurations are indicated as icons in green and red.}
\label{Regions}
\end{figure}
\end{center}

These regions are shown in Figure \ref{Regions}.  

\begin{conj}
The only positive mass planar equilateral 5-body central configurations in the Newtonian case are the regular pentagon and star with equal masses and the symmetric configurations discussed above.
\end{conj}

\section{Acknowledgements}

The authors would like to thank Manuele Santoprete for the suggestion to study this class of configuration.

\section{Competing Interests}
The authors have no competing interests to declare.

\bibliography{CelMechEtc}{}
\bibliographystyle{amsplain}

\end{document}